\documentclass[12pt]{sven}

\usepackage{theorem, epsfig, amsfonts, amsmath, fancybox}

\newcommand{\theo}[2]{\vspace{\baselineskip} 
\noindent\framebox[\textwidth]{\parbox{0.95\textwidth}{
\begin{#1} \label{#1} #2 \end{#1} } } \vspace{\baselineskip} }

\theoremstyle{plain}
\theorembodyfont{\rmfamily}

\newtheorem{PIPA-conv}{Theorem}[section]
\newtheorem{PIPA-fail}{Lemma}[section]

\begin{document}

\pagestyle{myheadings}

\markboth{Sven Leyffer}
         {The Penalty Interior-Point Method Fails to Converge}

\title{The Penalty Interior-Point Method Fails to 
	Converge\thanks{Preprint ANL/MCS-P1091-0903}}

\date{September 22, 2003}

\author{Sven Leyffer\thanks{Mathematics and Computer Science Division, 
			Argonne National Laboratory, Argonne, IL 60439
			({\tt leyffer@mcs.anl.gov}).}}

\maketitle

\begin{abstract}

	\noindent
	Equilibrium equations in the form of complementarity conditions
	often appear as constraints in optimization problems. Problems of
	this type are commonly referred to as 
	mathematical programs with complementarity constraints (MPCCs). 
	A popular method for solving MPCCs is the penalty interior-point 
	algorithm (PIPA). This paper presents a small example for which PIPA  
	converges to a nonstationary point, providing a counterexample to 
	the established theory. The reasons for this adverse behavior are 
	discussed.


	\noindent
	{\bf AMS-MSC2000:} 90C30, 90C33, 90C51, 49M37, 65K10.
\end{abstract}

\section{Introduction}
\label{intro}
\setcounter{equation}{0}

Equilibrium equations in the form of complementarity conditions
often appear as constraints in optimization problems. Problems of
this type are commonly referred to as 
mathematical programs with complementarity constraints
(MPCCs), and arise in many engineering and economic
applications; see the survey \cite{FerPan:97} and the monographs
\cite{LuoPanRal:96,OutKocZow:98} for further references. The growing 
collections of test problems \cite{MacMPEC,MPECWORLD} indicate
that this an important area. MPCCs can be expressed in general as
\be\label{MPEC}
\begin{array}{ll}
	\mini 	& f(x,y,w,z)			\\
	\st   	& x \in X			\\
		& F(x,y,w,z) = 0		\\
     		& 0 \leq y  \; \perp \;  w \geq 0,
\end{array}
\ee
where $X \subset \R^n$ is a polyhedral set, $f$ and $F$ are twice
continuously differentiable functions, and $y, w \in \R^m$. The 
complementarity constraint $y  \; \perp \;  w$ means that either 
a component of $y$ is zero, or the corresponding component of $w$ is 
zero, which implies that $y^T w = 0$.

Many approaches exist for solving MPCCs. These include branch-and-bound 
methods \cite{BarJF:88}, implicit nonsmooth approaches
\cite{OutKocZow:98}, piecewise SQP methods \cite{LuoPanRal:96}, 
perturbation and penalization approaches \cite{DirkFerrMeer:02}
analyzed in \cite{SchS:01}, and the penalty interior-point algorithm
(PIPA) \cite{LuoPanRal:96}, whose convergence properties are the subject 
of this note.

Recently, there has also been renewed interest in solving (\ref{MPEC})
with standard nonlinear programming (NLP) techniques by replacing the
complementarity constraint by 
\[ y^T w \leq 0, \; y \geq 0, \; w \geq 0 .\]
This approach had been regarded as numerically unsafe because (\ref{MPEC}) 
violates a constraint qualification at any feasible point \cite{SchSch:00}.
However, numerical \cite{FletLeyf:02} and theoretical evidence 
\cite{AnitM:00b,FLRS:02} suggests that sequential quadratic programming
(SQP) solvers are able to solve large classes of MPCCs. A similar 
analysis is being extended to interior-point methods (IPMs) 
\cite{BensShanVand:02,LiuSun:02,RaghBieg:02a}.

This note is organized as follows. The next section briefly describes
PIPA and presents its convergence result. Section~\ref{example} presents
a small counterexample  that shows that PIPA can fail to converge for 
certain admissible choices of parameters. Finally, Section~\ref{conclusion}
discusses the example and suggests a remedy that avoids this problem.

\section{The Penalty Interior-Point Algorithm (PIPA)}
\label{PIPA}

PIPA \cite[Chapter 6.1]{LuoPanRal:96} solves MPCCs of the form (\ref{MPEC}).
It has been studied in \cite{FukPan:98,Piep:01} 
and has been applied to a number of real applications in finance 
\cite{HuanPang:99}, target classification \cite{OlsoPangPrie:03}, 
and electricity markets \cite{HobbMetzPang:00}. 

PIPA generates a sequence of iterates $(x^k,y^k,w^k,z^k)$ that are 
strictly interior in $y^k > 0$, $w^k > 0$. The iterates are computed
by solving a quadratic direction-finding problem. In that sense, PIPA
combines aspects of interior-point and SQP methods.

We use $d$ to denote the step or displacement computed 
by PIPA. Subscripts such as $d_x$ refer the part of $d$ corresponding 
to the $x$-variables; superscripts are used to denote iterates or evaluation
of functions at a particular point, for example,
$\nabla f^k = \nabla f(x^k,y^k,w^k,z^k)$. Diagonal matrices are denoted 
by $W = \mbox{diag}(w)$ and $Y = \mbox{diag}(y)$ and $e = (1,\ldots,1)^T$. 

The algorithm solves the following direction-finding problem at every iteration:
\be\label{QP-dir}
	\begin{array}{ll}
	\mini 	& \nabla f^{k^T} d + \frac{1}{2} d_x^T Q^k d_x	\\
	\st   	& x^k + d_x \in X				\\
		& F^k + \nabla F^{k^T} d  = 0			\\
  		& Y^k d_w + W^k d_y 
		   = - Y^k w^k + \sigma \frac{y^{k^T}w^k}{m} e 	\\
		& \| d_x \|^2_2
		  \leq c \left(\| F^k \| + y^{k^T}w^k \right) ,
	\end{array}
\ee
where $c > 0$ and $\sigma \in (0,1)$ are parameters and $Q^k$ approximates
second-order information in $f$ and $F$. Note that (\ref{QP-dir}) has a 
quadratic trust-region-type constraint. If this constraint is replaced by
the $\ell_{\infty}$ norm, then (\ref{QP-dir}) becomes a quadratic program
(QP). Note also that the linearization of the complementarity constraint
$Yw = 0$ is relaxed by $\sigma \frac{y^{k^T}w^k}{m} e$, giving PIPA 
an interior-point flavor. Note that PIPA is not strictly an
interior-point algorithm, as it remains interior only with respect
to the complementary variables $y, w$ but not with respect to the other
constraints. We will show that the final trust-region type constraint is 
the source of the problems affecting PIPA. 

After determining a search direction, PIPA performs a backtracking line-search 
to ensure that $y^k, w^k$ remain positive and close to the central path.
This is achieved by finding the root $\tau$ of
\be\label{E-g} 
   g_k(\tau) = (1 - \rho) \sigma \frac{y^{k^T}w^k}{m} 
	     + \tau \left( \min_{1 \leq i \leq m} d^k_{y_i} d^k_{w_i}
			- \rho \frac{d_y^{k^T} d^k_w}{m} \right) .
\ee
Global convergence is enforced by further reducing the step length to
yield sufficient reduction in a quadratic penalty function
\be\label{E-P}
	P_{\alpha}(x,y,w,z) = f(x,y,w,z) 
			    + \alpha \left( \| F(x,y,w,z) \|^2
						    + y^T w \right),
\ee	
where the penalty parameter $\alpha>0$ is chosen to be the smallest $\alpha$
such that the step ensures sufficient model decrease,
\be\label{E-ModDecr}
\nabla f^{k^T} d - \alpha^p (1-\sigma) y^{k^T} w^k 
   \; < \; - \alpha^p (1-\sigma) y^{k^T} w^k
   \; < \; - y^{k^T} w^k.
\ee
Note that the penalty function $P_{\alpha}$ mixes a quadratic penalty for
the nonlinear equation $F$ with an $\ell_1$ penalty for complementarity
(since $y^k, w^k > 0$, it follows that $\| y^{k^T} w^k \|_1 = y^{k^T} w^k$).
Further details of the algorithm can be found in \cite{LuoPanRal:96}. PIPA
can be summarized as follows. 

\bigskip \noindent
{\bf Penalty Interior-Point Algorithm (PIPA)}				

\vspace{-1ex}
\noindent
\begin{tabbing} \hspace{0.5cm} \= \hspace{0.5cm} \= \hspace{0.5cm} \= \kill
1. {\em Initialization:\/}						\\
   \> Choose parameters $c > 0$, $\sigma, \gamma, \rho \in (0,1)$,		\\
   \> choose a starting point $(x^0,y^0,w^0,z^0)$ such that $y^0, w^0 > 0$
	suitably centered,							\\
   \> choose a penalty parameter $\alpha > 0$ and set $k = 0$.		\\
{\bf REPEAT}								\\
   \> 2. \> {\em Direction finding problem:\/}
            Solve problem (\ref{QP-dir}) for a trial step $d = d^k$.		\\
   \> 3. \> {\em Step size determination:\/} Find a step size $\tau = \tau_k$ 
									to:	\\
   \> \> 3.1 \> Ensure centrality and positivity of $(y,w)$ by finding 
		the root of $g_k(\tau)$ in (\ref{E-g})			\\
   \> \> \>  	or setting $\tau = 1 - \epsilon$ if this root does not exist
		or is greater than 1.					\\
   \> \> 3.2 \> Ensure sufficient reduction in the quadratic penalty function 
		$P_{\alpha}(x,y,w,z)$ 					\\
   \> \> \>	in (\ref{E-P}) by performing an Armijo-search on $P_{\alpha}$.\\
   \> \> Let $\tau_k$ be the step size determined in 3.1 and 3.2.	\\
   \> 4. \> {\em Update:\/}
	    $(x^{k+1},y^{k+1},w^{k+1},z^{k+1}) = (x^k,y^k,w^k,z^k) 
		 + \tau_k (d_x^k,d_y^k,d_w^k,d_z^k)$,  $k = k+1$.		\\
{\bf UNTIL} $\| d \| \leq \epsilon$
\end{tabbing}

\bigskip
This presentation of PIPA is less sophisticated than \cite{LuoPanRal:96},
which allows for instance $\sigma$ to vary with the iteration. 
Convergence of PIPA is established in \cite{LuoPanRal:96} under the following 
two assumptions.
\begin{description}
  \item{{\bf [SC]}} Strict complementarity of the solution, namely $y^* + w^* > 0$.
  \item{{\bf [NS]}} Nonsingularity of the matrix 
	\[ \left[ \begin{array}{ccc} \nabla_y F^* & \nabla_w F^* & \nabla_z F^* \\
				     W^*	  & Y^* 	 & 0	
		  \end{array} \right] .
	\]
\end{description}
A sufficient condition of {\bf [NS]} is that the Jacobian of $F$ satisfies
the mixed $P_0$ property. The following theorem summarizes the convergence
results in \cite{LuoPanRal:96}.

\theo{PIPA-conv}{{\bf \cite[Theorem 6.1.17]{LuoPanRal:96}}~~~Suppose that the 
		Hessian matrices $W^k$ are bounded and that $\sigma > 0$. If the 
		penalty parameter $\alpha$ is bounded, then every limit point of 
		$(x^k,y^k,w^k,z^k)$ that satisfies {\bf [SC]} and {\bf [NS]} 
		is a stationary point of (\ref{MPEC}).}

In the next section, we present a small example for which PIPA converges
to a non-stationary point, contradicting Theorem~\ref{PIPA-conv}.

\section{A Counterexample}
\label{example}
\setcounter{equation}{0}

In this section, we show that PIPA may fail to converge to a 
stationary point. This result contradicts the convergence results of 
\cite{LuoPanRal:96}. This section examines the behavior of 
PIPA applied to the following small example:
\be\label{Example}
\begin{array}{ll}
	\mini 	& x + w				\\
	\st   	& -1 \leq x \leq 1		\\
		& -1 + x + y = 0		\\
     		& 0 \leq y  \; \perp \;  w \geq 0.
\end{array}
\ee
The solution of (\ref{Example}) is $(x,y,w)^* = (-1,2,0)$, which
satisfies the assumptions {\bf [SC]} and {\bf [NS]}, since $y^* + w^* = 2 > 0$
and
\[ \left[ \begin{array}{ll} \nabla_y F^* & \nabla_w F^* \\
			    W^*		 & Y^*	
	  \end{array} \right] \; = \;
   \left[ \begin{array}{ll} 1 & 0 \\ 0 & 2 \end{array} \right] 
\]
is nonsingular. Next, we show that this example generates a limit point that
is not stationary. Specifically, we perform a standard iteration of the 
algorithm and bound the limit away from stationarity. The following lemma
summarizes our result.

\theo{PIPA-fail}{PIPA applied to the Example (\ref{Example}) and started at 
		 $(x^0,y^0,w^0) = (0,1,0.02)$, with $c  = 1, \sigma = 0.1$, 
		 $\gamma = 0.01, \rho = 0.9$ and $\alpha = 2$ generates a 
		 sequence of iterates satisfying 
		 \bea
			1 \leq y^k \leq y^{k+1}	&	& \label{E-ind-1}	\\
			w^{k+1} \leq \frac{1}{2} w^k \leq \frac{2}{100}
						&	& \label{E-ind-2}	\\
			x^{k+1} \geq x^k - \sqrt{y^k w^k} > -1	
						& 	& \label{E-ind-3}	\\
			y^{k+1}w^{k+1} \leq \frac{1}{2} y^k w^k	
						&	& \label{E-ind-4}	\\
			\frac{5}{9} \leq \tau_k \leq 1 .
						&	& \label{E-ind-5}	
		 \eea }

\noindent
{\bf Proof.} The starting point satisfies the linear equation,
which implies that they are satisfied for all subsequent iterations.
The theorem follows by induction. Clearly, the assertions hold for the 
starting point. Now assume that (\ref{E-ind-1}) to (\ref{E-ind-5}) hold for 
$k-1$, and show that they also hold for $k$.

The direction-finding problem for this example can be simplified as follows:
\be\label{Ex-QP} \begin{array}{ll}
	\mini 	& d_x + d_w					\\
	\st   	& -1 \leq x^k + d_x \leq 1			\\
		& d_x + d_y = 0					\\
     		& w^k d_y + y^k d_w = (\sigma - 1) y^k w^k	\\
		& - \sqrt{y^k w^k} \leq d_x \leq \sqrt{y^k w^k} ,
\end{array}
\ee
where $Q^k = 0$ for simplicity. However, the example remains valid for 
certain positive semi-definite bounded $Q^k$.

First we show that the optimal basis $B$ to the LP (\ref{Ex-QP}) is given by
\[ \left[ \begin{array}{ccc} 1 & 0 & 1 \\ 1 & w^k & 0 \\ 0 & y^k & 0 
	  \end{array} \right]^{-1}
   \; = \;
   \left[ \begin{array}{ccc} 	0 & 1 & -\frac{w^k}{y^k} \\
				0 & 0 & \frac{1}{y^k} 	\\
				1 & -1& \frac{w^k}{y^k} \end{array} \right] .
\]
From this, multipliers can be computed as
\[ \lambda^k \; = \; B^{-1} c \; = \; 
	\left[ \begin{array}{ccc} 0 & 1 & -\frac{w^k}{y^k} \\
				  0 & 0 & \frac{1}{y^k} 	\\
				  1 & -1& \frac{w^k}{y^k} \end{array} \right] \;
	\left( \begin{array}{c} 1 \\ 0 \\ 1 \end{array} \right) \; = \;
  	\left( \begin{array}{c}  -\frac{w^k}{y^k} \\ \frac{1}{y^k} \\ 
				 1 + \frac{w^k}{y^k} \end{array} \right).
\]
Optimality follows, since $\lambda^k_3 \geq 0$ for the only inequality
constraint. The step $d$ can also be computed as
\[ d^k \; = \; B^{-T} b^k
       \; = \; \left( \begin{array}{c} - \sqrt{y^k w^k} \\
					\sqrt{y^k w^k}	\\
			- \frac{9}{10} w^k - \frac{w^k}{y^k} \sqrt{y^k w^k}
		\end{array} \right) ,
\]
where $b^k = (0, -\frac{9}{10} y^k w^k , - \sqrt{y^k w^k})^T$ 
is the right-hand side of the active constraints. 

Note that for any positive semi-definite $Q^k$, this solution remains
unchanged as long as $\lambda^k_3 \geq 0$. The contribution of the 
Hessian to the multipliers can be estimated as 
\[ \tilde{\lambda_3^k} = \lambda_3^k + B^{-1} Q^k d^k 
			\geq 1 + \frac{w^k}{y^k} - \| B^{-1} Q^k \| \| d^k \|,
\]
which remains positive for $\| Q^k \|$  or $\| d^k \|$ sufficiently small. 
Thus, the same conclusions apply even if a positive definite $Q^k$ is used.

Next, we show that (\ref{E-ind-3}) follows by applying the step $d^k$, 
using $\tau^k \leq 1$ and induction:
\bea 
	x^{k+1} & =    & x^k + \tau^k d^k_x 				\nonumber \\
		& \geq & x^k - \sqrt{y^k w^k}				\nonumber \\
		& \geq & \dps x^k - \sum_{l=0}^k \left(\frac{1}{2}\right)^l 
					    \sqrt{y^0 w^0} 		\nonumber \\
        	& =    & \frac{- \sqrt{2}}{10} \left(\frac{1 - 
				\left(\frac{1}{\sqrt{2}}\right)^k}
			       {1 - \frac{1}{\sqrt{2}}} \right) .  	\label{E-xk}
\eea
It is easy to show that the expression in the last line is greater than $-1$,
which implies that the ``choice'' of active set in the LP (\ref{Ex-QP}) was
correct, and the lower bound on $x$ is never active during the iteration.
Next we show that (\ref{E-ind-1}) follows by applying the step in the $y$ component,
\[ y^{k+1} = y^k + \tau^k \sqrt{y^k w^k} \geq y^k \geq 1. \]
Next consider (\ref{E-ind-5}), and observe that $\tau^k \leq 1$ follows trivially.
To obtain the lower bound, observe that the root of $g_k(\tau)$ in (\ref{E-g}) is
given by the following expression, where the superscripts $k$ have been
omitted for the sake of simplicity
\[ \begin{array}{lcl}
     \tau	& = & \dps - \frac{1}{10} \frac{y w}{d_w d_y} 
		    \; = \; \frac{y^2}{\left(9 y + 10 \sqrt{y w}\right) 
		            \sqrt{y w}}					\\
		& = & \dps \frac{1}{9} \frac{\sqrt{y}}{\sqrt{w}}
		      \frac{1}{1+\frac{10}{9}\frac{\sqrt{w}}{\sqrt{y}}}	
		    \; = \; \frac{1}{9} \frac{\sqrt{y}}{\sqrt{w}}
			\sum_{l=0}^{\infty}\left(-\frac{10}{9}
				\frac{\sqrt{w}}{\sqrt{y}} \right)^l	.
   \end{array}
\]
From the fact that $w \leq \frac{2}{100}$, $y \geq 1$, induction, and
$\frac{100}{729}\frac{\sqrt{w}}{\sqrt{y}} \geq 0$,
the step size can be bounded below by its leading terms:
\[ \begin{array}{lcl}
   \tau & \geq & \frac{1}{9} \frac{1}{\sqrt{\frac{2}{100}}} - \frac{10}{81}
		 + \frac{100}{729}\frac{\sqrt{w}}{\sqrt{y}}
		 - \frac{1000}{6561}\frac{\frac{2}{100}}{1} + \ldots	\\
	& \geq & \frac{10}{9}\frac{1}{\sqrt{2}} - \frac{10}{81} - \frac{20}{6561} \\
	& \geq & \frac{5}{9}.
   \end{array}
\]
It remains to show that this lower bound on $\tau$ also passes the second
line-search criterion. The penalty update rule in \cite[(6.1.24)]{LuoPanRal:96} 
finds the smallest integer $p \geq 1$ such that
\[ \nabla f^{k^T} d - \alpha^p (1-\sigma) y^{k^T} w^k 
   \; < \; - \alpha^p (1-\sigma) y^{k^T} w^k
   \; < \; - y^{k^T} w^k 
\]
holds, where we have used the fact that $\|F^k\| \equiv 0$ for all $k$.
Since $\nabla f^{k^T} d < 0$ and $\alpha (1-\sigma) > 1$, this
is always satisfied for $p = 1$, and the penalty parameter is never increased.
Finally, this step size also satisfies the conditions
for the Armijo search in Step 3.2. The actual reduction clearly satisfies
\[ P_{\alpha}^{k+1} - P_{\alpha}^k 
   = \tau d_x + \tau d_w + \alpha \left( y^{k+1} w^{k+1} - y^k w^k \right).
\]
Since $y^{k+1} w^{k+1} \leq \frac{1}{2} y^k w^k$, it follows that
\[ P_{\alpha}^{k+1} - P_{\alpha}^k
   \leq \tau d_x + \tau d_w - \alpha \frac{1}{2} y^k w^k,
\]
which implies the sufficient reduction condition, 
\[ P_{\alpha}^{k+1} - P_{\alpha}^k \leq \gamma \tau d_x + \gamma \tau d_w 
					- \alpha (1 - \sigma) \gamma y^k w^k ,
\]
since $\gamma \leq 1$ and $(1 - \sigma) \gamma \leq \frac{1}{2}$. 

Next consider (\ref{E-ind-2}). Since $\tau \geq \frac{5}{9}$, it follows that
\be\label{E-wnew} 
   w^{k+1} \; = \; w^k 
		 + \tau \left( - \frac{9}{10} w^k - \frac{w^k}{y^k} \sqrt{y^k w^k}
			\right)
	   \; \leq \; \frac{1}{2} w^k - \frac{5}{9} \frac{w^k}{y^k} \sqrt{y^k w^k} .
\ee
Thus, it is possible to bound the new complementarity violation after the 
step by
\[ \begin{array}{lcl}
   y^{k+1} w^{k+1} & \leq & \left( y^k + \sqrt{y^k w^k} \right)
			    \left( \frac{1}{2} w^k - \frac{5}{9} 
				   \frac{w^k}{y^k} \sqrt{y^k w^k} \right)	\\
		   & = & \frac{1}{2} y^k w^k - \frac{5}{9} w^k \sqrt{y^k w^k}
			 + \frac{1}{2} w^k \sqrt{y^k w^k} - \frac{5}{9} (w^k)^2 \\
		   & \leq & \frac{1}{2} y^k w^k ,
   \end{array}
\]
obtaining (\ref{E-ind-4}). This concludes the proof.	\hfill $\Box$

\bigskip
By taking the limit ($k \rightarrow \infty$) in (\ref{E-xk}), it follows, 
that all control iterates satisfy the bound
\be \label{E-xlim}
    x^{k+1} \; \geq \; \frac{- 2}{10 \left(\sqrt{2} - 1\right)} \approx - 0.4828 .
\ee
Similarly, an upper bound can be derived for the state iterates as
\[
  y^{k+1} 
    \; \leq \; 
    y^0 + \sum_{l=0}^k \left(\frac{1}{\sqrt{2}}\right)^l \sqrt{y^0 w^0}
    \: = \: 
    1 + \frac{\sqrt{2}}{10} \left(\frac{1 - \left(\frac{1}{\sqrt{2}}\right)^k}
			       {1 - \frac{1}{\sqrt{2}}} \right), 
\]
which shows that
\be \label{E-ylim}
   y^{k+1} \; \leq \; 1 + \frac{2}{10 \left(\sqrt{2} - 1\right)} \approx 1.4828 .
\ee
These bounds can also be verified numerically. 
Example (\ref{Example}) and PIPA have been implemented in AMPL \cite{FouGayKer:93}
by using the ``looping extension'', which allows the convenient implementation
of algorithms in AMPL. The sequence that is generated is given in 
Table~\ref{T-Res} and confirms the bounds given above.
The columns headed ``ared'' and ``pred'' in Table~\ref{T-Res} refer to the actual
and predicted reduction. Their implications are discussed below.

\begin{table}[htb]
\begin{center}
\begin{tabular}{r|rrr|rr}
$k$ & $x^k$	   & $y^k$	  & $w^k$	   & ared	& pred		\\ \hline
1 &            0   &            1 &           0.02 &  		& 		\\ 
2 & -0.096022613   &    1.0960226 &   0.0058578644 &  -0.198      &  -0.137	\\ 
3 &  -0.17606958   &    1.1760696 &  0.00016323495 &  -0.0974     &  -0.0982	\\ 
4 &  -0.18991126   &    1.1899113 &  1.4549224E-05 & -0.0143     &  -0.0143	\\ 
5 &   -0.1940679   &    1.1940679 &  1.4171928E-06 & -0.00421    &  -0.0042	\\ 
6 &  -0.19536745   &    1.1953675 &  1.4145236E-07 & -0.00131    &  -0.0013	\\ 
7 &  -0.19577825   &    1.1957782 &  1.4223933E-08 & -0.000412   & -0.000411 	\\ 
8 &  -0.19590853   &    1.1959085 &   1.433645E-09 & -0.00013    & -0.00013	\\ 
9 &   -0.1959499   &    1.1959499 &  1.4460519E-10 & -4.14e-05   & -4.14e-05	\\ 
10 &  -0.19596304  &      1.195963 &  1.4586827E-11& -1.32e-05   & -1.31e-05
\end{tabular}
\caption{Iterates, actual and predicted reduction of AMPL implementation of PIPA}
\label{T-Res}
\end{center}
\end{table}

Thus, it follows that all iterates remain in a compact set, namely,
\[ \begin{array}{rcccccl}
	0 		& = & x^0 & \geq & x^k & \geq & -0.4828	\\
	1 		& = & y^0 & \leq & y^k & \leq &  1.4828	\\
	\frac{2}{100}	& = & w^0 & \geq & w^k & \geq &  0 ,
   \end{array}
\]
and the sequence $(x,y,w)^k$ has a limit point $(x,y,w)^{\infty}$. From
(\ref{E-wnew}) it follows that 
\[ w^{k+1} \; \leq \; w^k \left( \frac{1}{2} - \frac{5}{9 y^k} \sqrt{y^k w^k} \right)
           \; \leq \; \frac{1}{2} w^k , 
\]
so that $w^{\infty} = 0$ and $y^{\infty^T} w^{\infty} = 0$. Therefore, PIPA converges 
to a limit point $(x,y,w)^{\infty}$ such that
\[ \begin{array}{rcccccl}
	0 		& \geq & x^{\infty} & \geq & -0.4828	\\
	1 		& \leq & y^{\infty} & \leq &  1.4828	\\
	 		&      & w^{\infty} & =    &  0 .
   \end{array}
\] 
This limit point satisfies complementarity and the linear constraints
for problem (\ref{Example}) (since $(x,y,w)^0$ satisfies the linear constraints). 
Moreover, this limit point also satisfies the assumptions of 
Theorem~\ref{PIPA-conv}. Assumption {\bf [SC]} is satisfied because
$y^{\infty} + w^{\infty} \geq 1 > 0$, and Assumption {\bf [NS]} holds because
\[ \left[ \begin{array}{ll} M		& I		\\
			    W^{\infty}	& Y^{\infty}	
	  \end{array} \right] \; = \;
   \left[ \begin{array}{ll} 1 & 0 \\ 0 & y^{\infty} \end{array} \right] 
\]
is nonsingular for all $1 \leq y^{\infty} \leq 1.4828$. However, the limit point
is clearly not a stationary point, thereby contradicting Theorem~\ref{PIPA-conv}.

\section{Discussion and Remedy}
\label{conclusion}

In the preceding section, a small example has been presented for which PIPA 
converges to a feasible but nonstationary limit point. This example 
contradicts \cite[Theorem 6.1.17]{LuoPanRal:96}, which establishes 
the convergence of PIPA. The example might be considered contrived in the
sense that we choose a particular starting point and parameter values. 
However, for any different starting point it is possible to construct
parameters such that PIPA fails.

The reason for this apparent failure of PIPA is the trust-region type 
constraint
\be\label{E-TR} 
   \| d_x \|^2_2 \; \leq \; c \left(\| F^k \| + y^{k^T}w^k \right) 
               \; =: \; \Delta_k
\ee
in the direction-finding problem (\ref{QP-dir}). The trust-region radius 
$\Delta_k$ converges to zero as the iterates approach feasibility, thereby 
limiting the progress toward optimality in the controls $x$. This adverse 
behavior can be expected to occur whenever the iterates approach feasibility
``faster'' than optimality.

Linking the trust-region radius to the feasibility of the problem is 
counterintuitive. Normally, $\Delta_k$ is controlled by
the algorithm, rather than by the iterates. In particular, it should
reflect how well the model problem (\ref{QP-dir}) approximates the
original problem, measured, for instance, by the agreement between actual
reduction,
\[ \mbox{ared } \; = \; P_{\alpha}^k - P_{\alpha}^{k+1} , \]
and predicted reduction,
\[ \mbox{pred } \; = \; \nabla f^{k^T} d - \alpha^p (1-\sigma) y^{k^T} w^k . \]
Table~\ref{T-Res} shows that in this example, both the actual and the predicted 
reductions show near perfect agreement toward the end. Thus, one would expect 
$\Delta_k$ to be increased rather than decreased.
The adverse situation in which (\ref{E-TR}) limits progress toward
optimality can be easily detected. The Lagrange multiplier of (\ref{E-TR})
indicates whether progress can be made by relaxing this constraint.
In the present example, this multiplier converges to 1.

One possible way to remedy the adverse behavior of PIPA is to control
the trust-region radius by using standard trust-region techniques. The 
direction-finding problem then becomes
\be\label{QP-dir1}
	\begin{array}{ll}
	\mini 	& \nabla f^{k^T} d + \frac{1}{2} d_x^T Q^k d_x	\\
	\st   	& x^k + d_x \in X				\\
		& F^k + \nabla F^{k^T} d  = 0			\\
  		& Y^k d_w + W^k d_y 
		   = - Y^k w^k + \sigma \frac{y^{k^T}w^k}{m} e 	\\
		& \| d_x \| \leq \Delta_k,
	\end{array}
\ee
where the trust-region is now controlled by the algorithm. A typical updating 
scheme for $\Delta_k$ is as follows. First compute the ratio of actual over
predicted reduction
\[ \rho_k := \frac{\mbox{ared}}{\mbox{pred}} , \]
and then update $\Delta_k$ according to
\[
  \Delta_{k+1} \in \left\{ \begin{array}{ll}
		\left[\gamma_0 \Delta_k,\gamma_1 \Delta_k\right] 
			& \mbox{if } \rho_k < \eta_1  \\
		\left[\gamma_1 \Delta_k,\Delta_k\right] 
			& \mbox{if } \rho_k \in [\eta_1,\eta_2) \\
		\left[\Delta_k,\gamma_2 \Delta_k\right] 
			& \mbox{if } \rho_k \geq \eta_2
		\end{array} \right.
\]
in Step 4, where $0 < \gamma_0 \leq \gamma_1 < 1 \leq \gamma_2$, and 
$0 < \eta_1 < \eta_2 < 1$. Typical values for these constants are 
$\gamma_0 = \gamma_1 = 0.5$, $\gamma_2 = 2$, $\eta_1 = 0.25$, 
and $\eta_2 = 0.75$.

This results in a mixed trust-region/line-search algorithm. Convergence
can be established under the assumption that the Jacobian of $F$,
\[ \left[ \begin{array}{ccc} \nabla_y F & \nabla_w F & \nabla_z F \\
			     W	  & Y 	 & 0	
	  \end{array} \right] ,
\]
remains nonsingular and has a bounded inverse along the lines of the 
coupled trust-region approach of Dennis et al. \cite{DennHeinVice:98},
since the QP subproblem can be regarded as a problem in $d_x$ only after
eliminating $(d_y, d_w, d_z)$ using the linearization of $F$.

\section*{Acknowledgments}

I am grateful to Roger Fletcher and Jorge Nocedal for many fruitful 
discussions on MPCCs and interior-point methods. I am also grateful
to Jong-Shi Pang and an anonymous referee for useful comments 
on an earlier version of the paper.
	This work was supported by the Mathematical, Information, and
	Computational Sciences Division subprogram of the Office of
	Advanced Scientific Computing Research, Office of Science,
	U.S. Department of Energy, under Contract W-31-109-ENG-38.


\vfill
\begin{flushright}
\scriptsize
\framebox{\parbox{2.4in}{The submitted manuscript has been created
by the University of Chicago as Operator of Argonne
National Laboratory ("Argonne") under Contract No.\
W-31-109-ENG-38 with the U.S. Department of Energy.
The U.S. Government retains for itself, and others
acting on its behalf, a paid-up, nonexclusive, irrevocable
worldwide license in said article to reproduce,
prepare derivative works, distribute copies to the
public, and perform publicly and display publicly, by or on
behalf of the Government.}}
\normalsize
\end{flushright}

\end{document}